\pdfoutput=1

\documentclass[11pt]{article}
\usepackage{amssymb}
\usepackage{graphicx}

\parskip \baselineskip

\newcommand{\ol}{\overline}

\begin{document}
\title{Finding or counting all shellings of a simplicial complex}

\author{ Marcel Wild}

\date{}
%\begin{document}
\maketitle

\begin{quote}
A{\scriptsize BSTRACT}: {\footnotesize The shellability status of previously investigated simplicial complexes with up to 24 facets is settled. In case of shellability the exact number of shellings is determined. Our algorithm merely relies on the facets, and not on additional information such as the face-numbers.}
\end{quote}

\section{Introduction}

A survey on the relevance of shellings in algebra and combinatorics can be found in [B]. A {\it simplicial complex}, or henceforth just {\it complex} [B,p.228] is a family $\Delta \subseteq 2^W$ which is hereditary in the sense that for all $X \in \Delta$ it follows from $Y \subseteq X$ that $Y \in \Delta$. The members of $\Delta$ are called {\it faces} and the numbers $f_i$ of $i$-element faces are the {\it face-numbers} of $\Delta$.
The maximal members of $\Delta$ are its {\it facets}. One verifies that the intersection of complexes is again a complex. For any family ${\cal F}\subseteq 2^W$ the {\it generated}  complex ${\cal F}\!\!\downarrow$ is the family of all subsets of sets from ${\cal F}$. It will be handy to write $X\prec Y$ if $X$ is a subset of $Y$ of cardinality $|Y|-1$. A {\it partial shelling} of $\Delta$ is a sequence $(G_1, G_2, \cdots, G_m)$ of distinct facets such that for all $2 \leq k \leq m$ one has [B,p.229]:

(1) \quad $(\forall i<k)(\exists j<k)\ G_i\cap G_k\subseteq G_j\cap G_k\prec G_k$

 Thus (1) demands that the facets of the intersection complex $(\{G_1,\dots,G_{k-1}\}\!\!\downarrow)\cap(\{G_k\}\!\!\downarrow)$ all have cardinality $|G_k|-1$. If $m$ is the total number $n$ of facets then we speak of a {\it shelling}.
Checking the shellability of $\Delta$ by testing all $n!$ permutations of its $n$ facets is only feasible for very small values of $n$. The only nontrivial published alternative seems to be [MH] which relies on knowing the face numbers\footnote{Actually, what needs to be known are certain numbers $h_1,h_2,...$ which however, presupposing shellability, are as hard to get as the face numbers [B,p.231,(7.5),(7.6)].}.  Unfortunately calculating the face numbers of a  complex is NP-hard, and getting from the face numbers to a shelling as proposed in [MH] is not straightforward either. In a way it rather works the other way around: If we {\it have} a shelling of $\Delta$, then the face numbers come as a perk [B,p.231,(7.5),(7.7)]. In the present article we compute all (if any) shellings  of a complex from a mere\footnote{Notwithstanding mentioned perk, our method seldom yields the face numbers of $\Delta$ because a {\it random} $\Delta$  is {\it unlikely} shellable. An altogether different approach to find the face-numbers from the facets can be found in [W4].} knowledge of its facets.

 Here comes the Section break-up. The abstract ideas in Section 2 about handling hereditary simple languages (e.g. greedoid languages) come to life in Section 3 where we replace the $n!$ permutations of $[n]:=\{1,2,\ldots,n\}$ (=index set of the facets) by certain pairs of type $(A,k)$ where $A\subseteq 2^{[n]}$ and $k\in [n]$. As motivated later, such an $(A,k)$ will be called a PSS. There are at most $n2^n\ (<<n!)$ many PSSes. Instead of storing them individually, compression techniques using wildcards will be applied.
For $\Delta$ to have a shelling it is necessary  that $|A|=n-1$ for at least one PSS $(A,k)$. In the latter case all shellings can in principle be enumerated one-by-one by depth-first searching the poset ${\cal P}$ of all PSSes. If there are too many shellings for enumeration then counting them may still be feasible. 

Our algorithms are applied to determine the precise number of shellings of some matroid complexes. (Matroid complexes are long known to be shellable.) As a sneak preview the matroid complex $\Delta$ of all trees of the complete graph $K_4$ has 16 facets and exactly 722965625856 shellings; the previously known lower bound  was 6!=720. Related to matroid complexes are chessboard complexes, but they need not be shellable, and indeed we prove the unshellability in some open cases.
Section 4 shows that, apart from shellings,  other kinds of restricted permutations can be counted or enumerated using appropriate PSS-posets. This e.g. concerns what we call 'peelings'. They are defined by a baby version of (1) and they generalize linear extensions of posets.

The paper in front of you is a preliminary draft which needs further trimming. Comments of algebraic combinatorists (=potential co-authors) are welcome. In particular proposals of concrete (i.e. with spelled out facets) moderate-size complexes with unknown shellability status.

\section{A tool for handling simple hereditary languages}

Subsection 2.1 introduces basic terminology, in 2.2 we define PSS-posets ${\cal P}$, and 2.3 to 2.5 show how the enumeration or counting of simple hereditary languages can be handled by means of ${\cal P}$.

{\bf 2.1} Following the terminology of [BZ,p.288] for a set $E$ denote by $E^*$ the free monoid of all words $\alpha$ over the alphabeth $E$. We denote the empty word by $\square$. The {\it length} $|\alpha|$ is the number of letters in $\alpha$, and the {\it support} $\tilde{\alpha}$ is the set of letters in $\alpha$. Thus if $\alpha=dbc$ and $\beta=cbdc$ then $|\alpha|\neq|\beta|$ but $\tilde{\alpha}=\tilde{\beta}=\{b,c,d\}$. A word $\alpha$ is {\it full} if $\tilde{\alpha}=E$, and it is  {\it simple} if it does not contain any letter more than once, i.e. $|\alpha|=|\tilde{\alpha}|$.
A {\it language} over $E$ is is a non-empty set ${\cal L}\subseteq E^*$; it is called {\it simple} if every word in ${\cal L}$ is simple. A language ${\cal L}$ is {\it hereditary} if from $\alpha=\beta b\in {\cal L}$ follows $\beta \in {\cal L}$
Each hereditary language ${\cal L}$ is determined by its {\it basic}, i.e. non-extendible words (which can be of different lengths). For instance, the language ${\cal L}_1\subseteq [4]^*$ defined by
the six basic words

(2)\quad $ 3412,\ 1432,\ 1342,\ 3142,\ 1423,\ 3421$

contains 6+6+4+2+1=19 words, counted according to decreasing length. We shall return to ${\cal L}_1$ later on. 

{\it Convention:} All considered languages are henceforth assumed to be simple and hereditary. Further we reserve $n$ for 
the cardinality of the alphabeth $E$, which often is $E=[n]$.

An {\it accessible} set system on a set $E$ is a non-empty subset ${\cal S}$ of the powerset $2^E$ such that for each $X\in {\cal S}$ there is $x\in X$ with $X\setminus\{x\}\in{\cal S}$. Then  ${\cal S}$ {\it induces} a language  ${\cal L(S)}$ as follows. For each (simple) word $a_1a_2\ldots a_k\in E^*$ it holds that

(3)\quad $a_1a_2\ldots a_k\in {\cal L(S)}\ :\Leftrightarrow\ (\forall 1\le i\le k)\ \{a_1,\ldots,a_i\}\in {\cal S}$

 It follows from (3) that the support of each word in ${\cal L(S)}$ is in ${\cal S}$, and each set in ${\cal S}$
occurs this way. Hence $\widetilde{\cal L(S)}={\cal S}$. Consequently each language ${\cal L}$ admits {\it at most  one} accessible set system  ${\cal S}$ inducing
${\cal L}$, namely ${\cal S}:=\tilde{\cal L}$.  We call a language 
{\it set-induced} if this happens. Hence there are as many set-induced languages with alphabeth $E$ as there are accessible set systems based on $E$, which is a tiny fraction of all simple hereditary languages. To fix ideas, let
 ${\cal L}_2$ be the language on $E=[3]$ whose basic words are $123$ and $312$. Then ${\cal L}_2$ is not set-induced because if it was it would coincide with
${\cal L(S)}$ for ${\cal S}=\tilde{\cal L}_2=\{\emptyset,\{1\},\{3\},\{1,2\},\{1,3\},\{1,2,3\}\}$. However, $132\in {\cal L(S)}\setminus {\cal L}_2$. (Exercise: Is ${\cal L}_1$ set-induced?)  Likely the most thoroughly investigated set-induced languages are {\it greedoid languages} [BZ,p.289], which will reoccur in 4.3.

{\bf 2.2} Let ${\cal L}$ be a language (set-induced or not) for which only an implicite description is available. Our task is to calculate a certain auxiliary 
poset $({\cal P},<)$ that facilitates the handling of ${\cal L}$, foremost the one-by-one enumeration or counting of its full words. (Our main, but not exclusive, concern will be languages whose full words match the shellings of a simplicial complex.) As to defining ${\cal P}$, consider a family ${\cal P}\subseteq 2^E\times E$ of pairs $(A,b)$ that satisfies these conditions:

\begin{description}
\item[(4a)] {\bf Scope:} For all $\alpha b\in {\cal L}$ it holds that $(\tilde{\alpha},b)\in {\cal P}$;
\item[(4b)] {\bf Normalization:} If $(\emptyset,b)\in {\cal P}$ then $ b\in {\cal L}$;
\item[(4c)] {\bf Linkage:} If $\alpha\in {\cal L}$ and  $(\tilde{\alpha},b),\ (\widetilde{\alpha b},c)\in {\cal P}$, then $\alpha b\in {\cal L}$
\end{description}

If (4a) to (4c) are fulfilled then ${\cal P}$, partially ordered by

(5) \quad $(A, b) < (A',b') : \ \Leftrightarrow \ A  \subsetneqq A'$,

is called a {\it PSS-poset} for ${\cal L}$, the acronym being explained in a moment. Observe that each language ${\cal L}$ admits PSS-posets, the most natural (but often illusive) one being \footnote{Since each word in ${\cal L}\setminus \{\square\}$ can be written as $\alpha b$, possibly with $\alpha=\square$.} 

 ${\cal P}_{\cal L}:=\{(\tilde{\alpha},b):\ \alpha b\in {\cal L}\}$. 

 Here we call $\tilde{\alpha}$ a {\it setment} because it stems from the {\it segment} $\alpha$, but $\tilde{\alpha}$ itself is a {\it set}. Furthermore, each pair of type $(\tilde{\alpha}, b)$ we name {\it setment-suffix}.  Let ${\cal P}$ be a PSS-poset for ${\cal L}$. Since ${\cal P}_{\cal L}\subseteq {\cal P}$ by (4a), some members of  ${\cal P}$ are setment-suffixes, but all $(A,b)\in {\cal P}$ will be referred to as {\it potential setment-suffixes (PSS)}, the first component $A$ being a {\it potential setment} (of $b$). For $s=0,1,\ldots,n-1$ the $s$-{\it level} of the PSS-poset consists of all PSSes $(A,k)$ with $|A|=s$.

\includegraphics[scale=0.62]{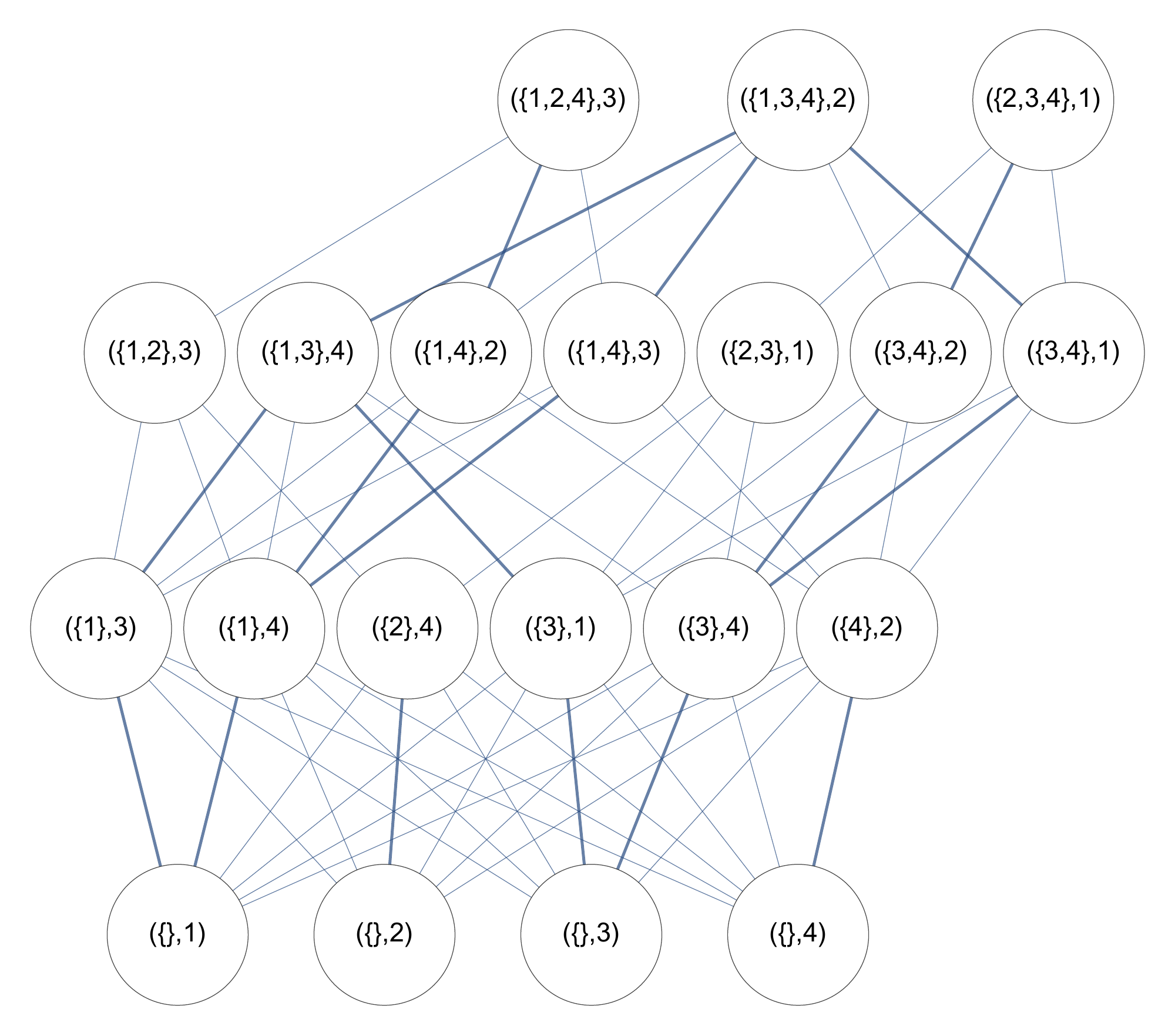}

{\it Figure 1: A PSS-poset for the language ${\cal L}_1$ defined by (2)}

In any PSS-poset $({\cal P},<)$ a {\it (rising) admissible} chain is an unrefinable chain of type

(6) \quad $(\{ \},a_1)<(\{a_1\},a_2)<\cdots<(\{a_1,\ldots,a_{k-2}\},a_{k-1})<(\{a_1,\ldots,a_{k-1}\},a_k)\hfill (k\ge 1)$. 

For instance, the chain $(\{\},3)<(\{3\},4)<(\{3,4\},1)$ in Figure 1 is admissible, but not the chain $(\{\},4)<(\{3\},4)<(\{3,4\},2)<(\{2,3,4\},1)$.
Dually a {\it falling admissible} chain is defined, such as $(\{1,2,4\},3)>(\{1,4\},2)$.

{\bf Proposition 1:} Let $({\cal P},<)$ be a PSS-poset for the language ${\cal L}$, i.e. (4a) to (4c) hold. Then the non-empty words in ${\cal L}$
 bijectively match the rising admissible chains in $({\cal P},<)$.

{\it Proof.} Consider first an admissible chain in $({\cal P},<)$, i.e. of type (6). By (4b) we have $a_1\in {\cal L}$. From $a_1\in {\cal L}$ and $(\tilde{a}_1,a_2)\in{\cal P}$ and
$(\widetilde{a_1 a_2},a_3)\in {\cal P}$ and (4c) follows
$a_1 a_2\in {\cal L}$, and so
forth until  $a_1\ldots a_k\in {\cal L}$. Conversely, start with any word $a_1\ldots a_k\in {\cal L}$. Then $(\{a_1,\ldots,a_{k-1}\},a_k)\in {\cal P}$ by (4a). Since ${\cal L}$ is hereditary, $a_1\ldots a_{k-1}\in {\cal L}$, hence again $(\{a_1,\ldots,a_{k-2}\},a_{k-1})\in {\cal P}$ by (4a), and so forth. This yields the claimed
 admissible chain.\hfill $\blacksquare$

By Proposition 1 the maximal admissible chains match the basic words of ${\cal L}$ which, recall, uniquely determine ${\cal L}$.
Often, as in Section 3, only the full (as opposed to basic) words of a language are of interest, and so we limit the upcoming discussion to this case. As we
show in 2.3, when the aim is {\it enumeration} of all full words, then depth-first search of the PSS-poset is best. If the full words merely need to be {\it counted}, then breadth-first search is better (2.4). If no full words exist then both types of search perform similarly to detect this (2.5). 

{\bf 2.3} By Proposition 1 enumerating all full words of ${\cal L}$, when a PSS-poset ${\cal P}$ is known, amounts to enumerate all admissible $n$-chains in 
${\cal P}$. The latter can be achieved by depth-first searching ${\cal P}$ either from above, or from below. Conceptually both barely differ but the algorithmic details {\it can} differ (as discussed later). The whole point of PSS-posets is to make implicite descriptions of languages more explicite. In this light, imagine that ${\cal L}_1$ (see (2)) was some implicitely rendered language whose basic words are unknown but for which the PSS-poset ${\cal P}_1$ in Figure 1 was somehow calculated. We now show how depth-first searching ${\cal P}_1$ from above leads to an explicite enumeration of the full words of ${\cal L}_1$.

\includegraphics[scale=0.87]{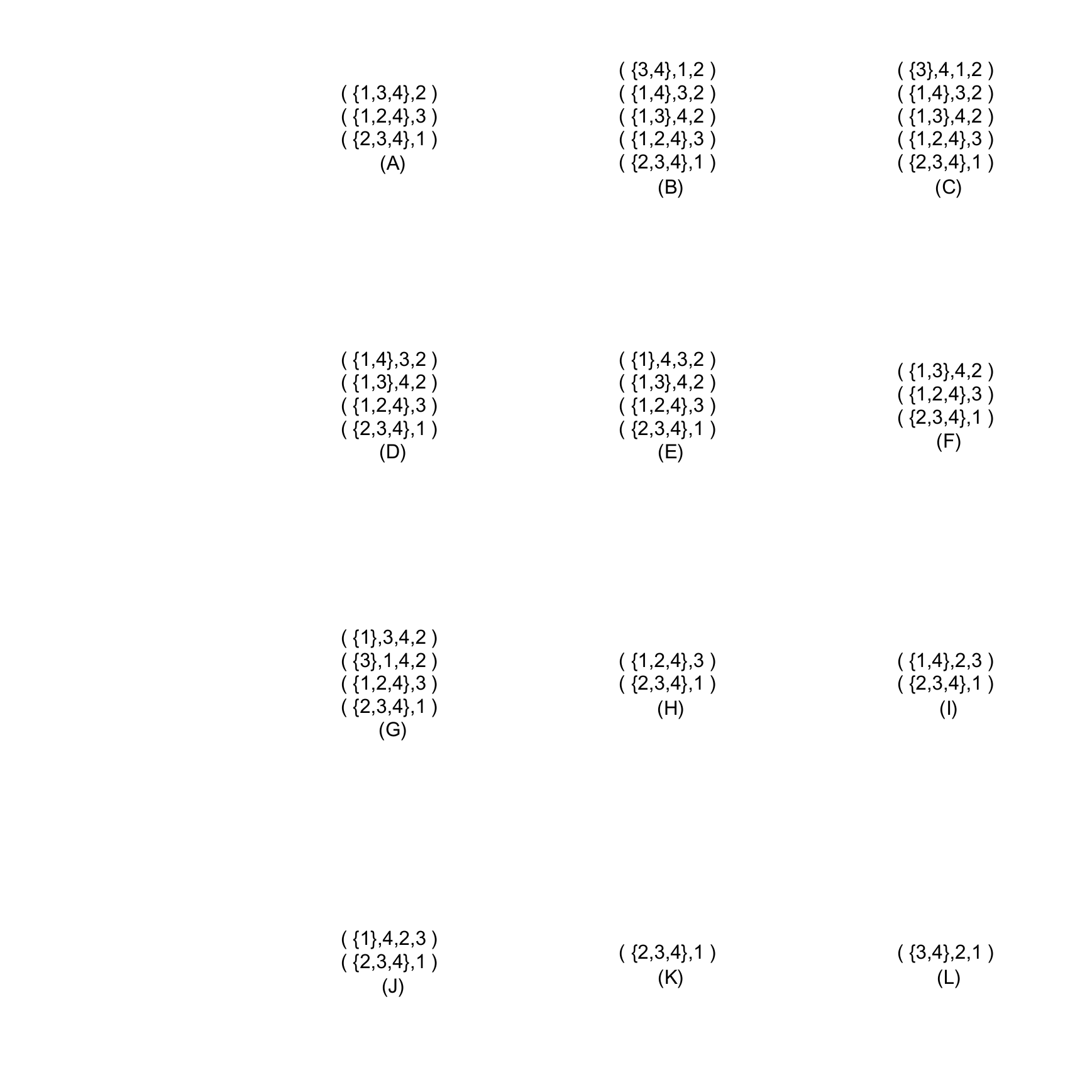}

{\it Figure 2: Depth-first search of the PSS-poset in Figure 1}

At the beginning (see (A)) the stack contains, in any order, all maximal element of the PSS-poset. Pick the top item $(\{1,3,4\},2)$ and determine those lower covers in the PSS-poset whose suffixes belong to the potential setment $\{1,3,4\}$. Thus $(\{1,3\},4)$ and $(\{1,4\},3)$ and
$(\{3,4\},1)$ qualify (indicated by thick lines), but neither $(\{1,4\},2)$ nor $(\{3,4\},2)$. Now replace  $(\{1,3,4\},2)$ by the qualifying lower covers which, however, get rewritten as
 $(\{1,3\},4,{\bf 2})$ and $(\{1,4\},3,{\bf 2})$ and $(\{3,4\},1,{\bf 2})$. In this way we remember that they originate from $(\{1,3,4\},{\bf 2})$; see (B). Again the order in which we put the qualifying lower covers on the stack is irrelevant. The new top item $(\{3,4\},1,2)$ is essentially treated the same way. That is, momentarily viewing it as $(\{3,4\},1)$, we pick those lower covers in the PSS-poset whose suffixes belong to $\{3,4\}$. Only $(\{3\},4)$ qualifies and this PSS, fleshed out as $(\{3\},4,1,2)$, kicks its father $(\{3,4\},1,2)$ from the stack; see (C). As soon as the potential setment of a top item has shrunk to a singleton, it represents a shelling (i.e. 3412 here) which is output; see (D). It is now clear how the stack keeps developing until it is empty after stage $(L)$. Summarizing we see that
the full words output at stages (C), (E),(G) (two), (J), and after (L), are exactly the ones in (2), as predicted.

\includegraphics[scale=0.67]{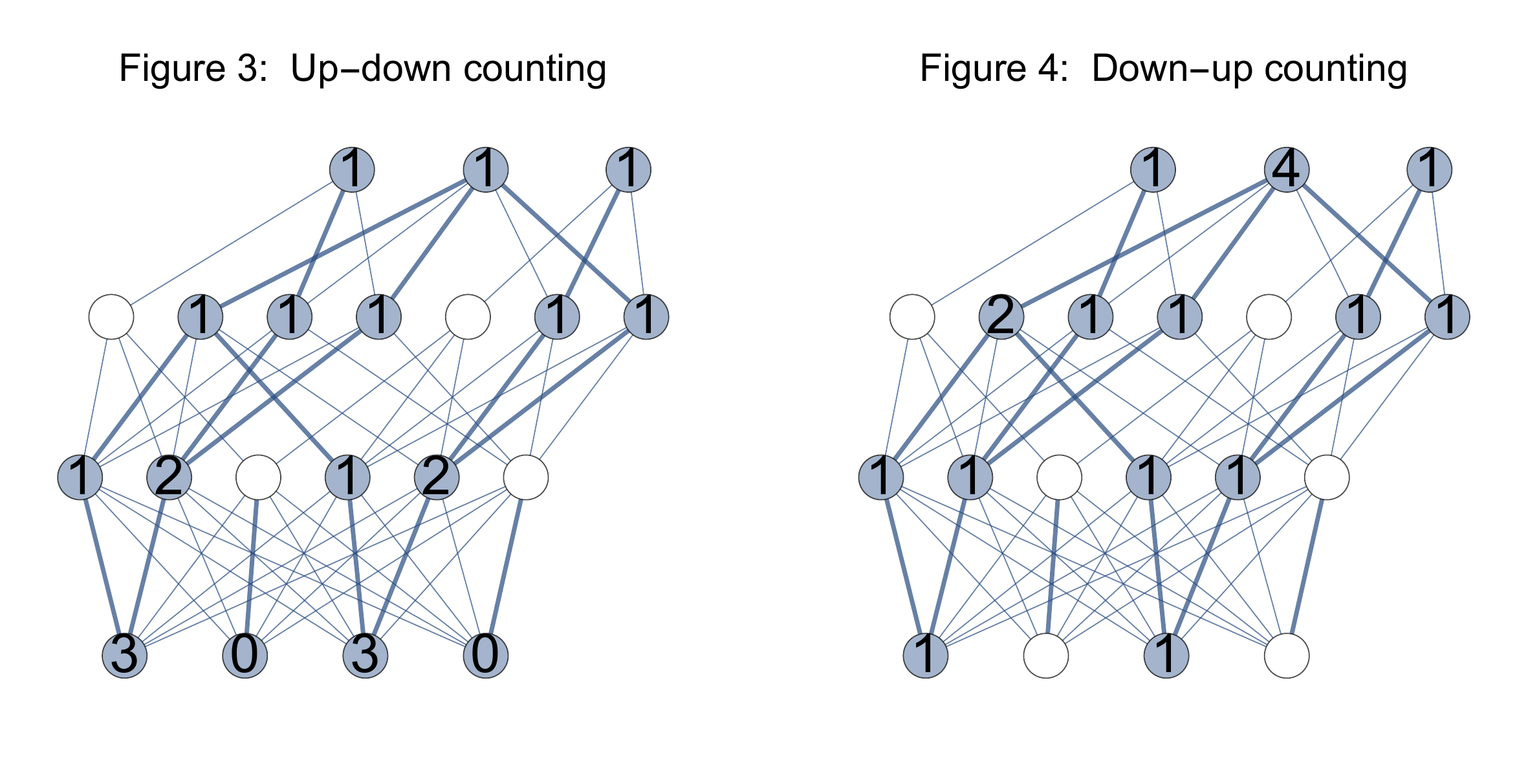}

{\bf 2.4}  As to counting full words (Figure 3), call a member of the 0-level {\it promising} if it has a {\it good} upper cover in the sense that the two are connected by a thick line. (The meaning of 'thick' lines was defined in 2.3.) For $s=0,1,\ldots, n-2$ by definition the {\it promising} setment-suffixes on the (s+1)-level (shaded grey, ignore the numbers) are the good upper covers $(A,b)$ of the promising s-level setment-suffixes, {\it unless} $(A,b)$ has no good upper covers itself. (That's why the third and sixth element of level 1 fail to qualify.)
For each promising $(A,b)$, except those on the top level, we record the set $guc(A,b)$ of its good upper covers. In Figure 3 most elements are promising, in contrast to the upcoming real-life examples. To summarize, using breadth-first search we have worked our way {\it up} in Figure 3. 

The way {\it down} begins with labeling all promising top level elements with the number 1. Then for $s=n-2,\ n-3,\ldots,0$ do the following. Label each promising s-level element $(A,b)$
with the sum of the labels of the elements in $guc(A,b)$. By inducting on $s$ the label of $(A,b)$ equals the number of falling admissible chains with last element $(A,b)$.
Hence in the end the sum of the labels on the 0-level equals the number of maximum admissible chains, i.e. by Proposition 1 the number of full words. Indeed 3+0+3+0=6. Specifically, three full words start with letter 1 and three with letter 3 (compare with (2)). 

Instead of going up-down, in  dual fashion one can do the counting by going down-up, see Figure 4. Thus we first go down, calculating good lower covers\footnote{There is an asymmetry between good upper and good lower covers of PSSes. All good upper covers of $(A,k)$ are of type $(A\cup\{k\},k')$, i.e. have the {\it same} setment $A\cup\{k\}$, whereas the good lower covers look like $A\setminus\{a_1\},\  A\setminus\{a_2\}$, i.e. have {\it different} setments. Accordingly it is likely more efficient (and we do so in Section 3) to target good {\it upper} covers.}, and then do the counting on the way up. One gets again 6 but this time as 1+4+1. Specifically, one and four and one full words end with the letters 3,2,1 respectively.

{\bf 2.5} What about merely deciding the existence of full words? Then either of breadth-first and depth-first search can be better. Specifically, if there are {\it no} full words (or at least this is suspected) then breadth-first search is better because it traces each thick line at most once\footnote{Using breadth-first search from below to trace all lines contained in {\it rising} admissible chains (say there are $N_1$ such thick lines) requires work proportional to $N_1$. Similarly if we use breadth-first search from above to trace the $N_2$ many thick lines featuring in {\it falling} admissible chains. Unfortunately, either of $N_1$ and $N_2$ may be much smaller than the other, and predicting which is next to impossible.}, whereas depth-first search traces each thick line at least once.
But if there {\it are} are full words (=admissible $n$-chains), in particular if they constitute a sizeable fraction of all $n$-chains in the PSS-poset, then depth-first is advised.

\section{The PSS-poset for the shelling problem}

In 3.1 an auxiliary algorithm that renders the models of a dual Horn formula in a compact format is presented. In Subsetion 3.2 we determine ad hoc all shellings of a toy complex $\Delta_0$. Subsections 3.3 to 3.5   draw on 3.1, 3.2 and the generalities of PSS-posets (Section 2) to taylor a PSS-poset for shellings of simplicial complexes.

{\bf 3.1}  We assume that the reader is familiar with basic Boolean logic. Thus recall that a clause is a disjunction of finitely many literals. One calls it a {\it Horn} clause if at most one literal is positive, and a {\it dual Horn} clause if at most one literal is negative. A conjunction of dual Horn clauses, such as

(7)\quad $\varphi(x_1,\ldots,x_8)\ =\ {\ol x}_1\wedge{\ol x}_2\wedge (x_5\vee {\ol x}_8)\wedge (x_3\vee {\ol x}_4\vee  x_5\vee x_7)$

is called a {\it dual Horn formula}, and it is these types of Boolean formulas that concern us here. The {\it model set} $Mod(\varphi)$ consists of all length 8 bitstrings $a=(a_1,\ldots,a_8)\in\{0,1\}^8$ with $\varphi(a)=1$. Let us see how $Mod(\varphi)$ can be calculated in compressed format in polynomial total time.

A $012e${\it -row} $r$ on a finite set $E$ is a quadruplet $r: = \{zeros(r), ones(r), twos(r), ebubbles(r)\}$ such that $E$ is a disjoint union of these four (possibly empty) sets. Further, if $ebubbles \neq \emptyset$ then it is a disjoint union of $t \geq 1$ many sets $eb_i$ (called $e$-$bubbles$) such that $\varepsilon_i: = |eb_i| \geq 2$ for all $1 \leq i \leq t$. Hence up to permutation of the entries $r$ can be visualized as

(8) \quad $r = \left(\underbrace{0, \cdots, 0}_{\alpha}, \underbrace{1, \cdots, 1}_{\beta}, \underbrace{2, \cdots, 2}_{\gamma}, \underbrace{e_1, \cdots, e_1}_{\varepsilon_1}, \cdots, \underbrace{e_t, \cdots, e_t}_{\varepsilon_t}\right)$

By definition $r$ {\it represents} (and will henceforth be identified with) the family of sets $A \subseteq W$ satisfying
$$A \cap zeros(r) = \emptyset \quad \mbox{and} \quad ones(r) \subseteq A \quad \mbox{and} \quad (\forall 1 \leq i \leq t) \ \ eb_i \cap A \neq \emptyset.$$
Put another way, if we identify  subsets $A \subseteq E$ with  bitstrings $a \in \{0,1\}^E$ in the usual way then $r$ contains exactly the following bitstrings $a:$ The $0$'s in $r$ are $0$'s in $a$, the $1$'s in $r$ are $1$'s in $a$, the {\it don't care} symbols $2$ in $r$ allow the corresponding bits in $a$ to be $0$ or $1$, and each wildcard $e_i e_i \cdots e_i$ in $r$ demands $a$ to have  ``at least one $1$ here''. Evidently $|r|  = 2^\gamma (2^{\varepsilon_1} -1) \cdots (2^{\varepsilon_t}-1)$. 

The $e$-algorithm\footnote{More precisely, the algorithm in [W1] is called {\it (implication) n-algorithm}. Strictly speaking in the present article we apply the {\it (dual implication) e-algorithm} since we are dealing with {\it dual} Horn formulas. Yet either algorithm can be replaced by the other, if afterwards the output is 'dualized' in the sense of switching in each row 0's and 1's, as well as $e$-bubbles and $n$-bubbles. (Unsurprisingly, $n...n$ means 'at least one 0 here'.)} from [W1] computes $Mod(\varphi)$  as a disjoint union of the three $012e$-rows in Table 1.

%From a technical (as opposed to conceptual) point of view the dual implication $e$-algorithm is similar to the {\it transversal e-algorithm} of [?], which delivers a compressed representation for the family of all transversals of a set system. See also 4.2.} from [?] 

\begin{tabular}{l|c|c|c|c|c|c|c|c|}
& $x_1$ & $x_2$ & $x_3$ & $x_4$ &  $x_5$ & $x_6$ &  $x_7$ & $x_8$ \\ \hline
 &  &  &  &  &   &   &   &   \\ \hline
$r_1=$ & 0 & 0 & 2 & 2 & 1 & 2 & 2 & 2   \\ \hline 
 $r_2 =$ & 0 & 0 & 0  & 0 & 0 & 2 & 0 & 0  \\ \hline
 $r_3=$ & 0 & 0 & $e$ & 2 & 0 & 2 & $e$ & 0 \\ \hline 
\end{tabular}

{\it Table 1: Compressed representation of $Mod(\varphi)$ for the dual Horn formula $\varphi$ in (7)}

With hindsight this is easy to verify. Namely, since $x_6$ doesn't appear in (7), we have $x_6=2$ in all rows. The 0's under $x_1,\ x_2$ are evident as well. Furthermore, if $x_5=1$ (this case is treated in $r_1$) then (7) is satisfied independent of the other variable values, whence the five don't care 2's. For $x_5=0$ there are two subcases, i.e. $r_2$ and $r_3$. In $r_2$ we have $x_3=x_7=0$, and so $x_4=0$. In $r_3$ we have $x_3=1$ or $x_7=1$ (neatly expressed by $ee$), and so $x_4=2$. Finally in both $r_2$ and $r_3$ the clause $x_5\vee{\ol x}_8$ forces $x_8=0$. The three $012e$-rows achieve a handy compression of altogether 32+2+12=46 models of $\varphi$.

{\bf 3.2} Here comes our first concrete simplicial complex. For $W=\{a,b,c,d,e,f,g\}$ consider  $\Delta_0\subseteq 2^W$ with facets

(9)\quad $F_1=\{a,c,d,f\},\ F_2=\{a,b,c,f\},\ F_3=\{a,b,c,d,e,g\},\ F_4=\{c,d,e,f,g\}$

Using shorthand notation (say $cf$ for $\{c,f\}$), here is the table of all intersections $F_i\cap F_j\ (i\neq j)$:

\includegraphics[scale=0.7]{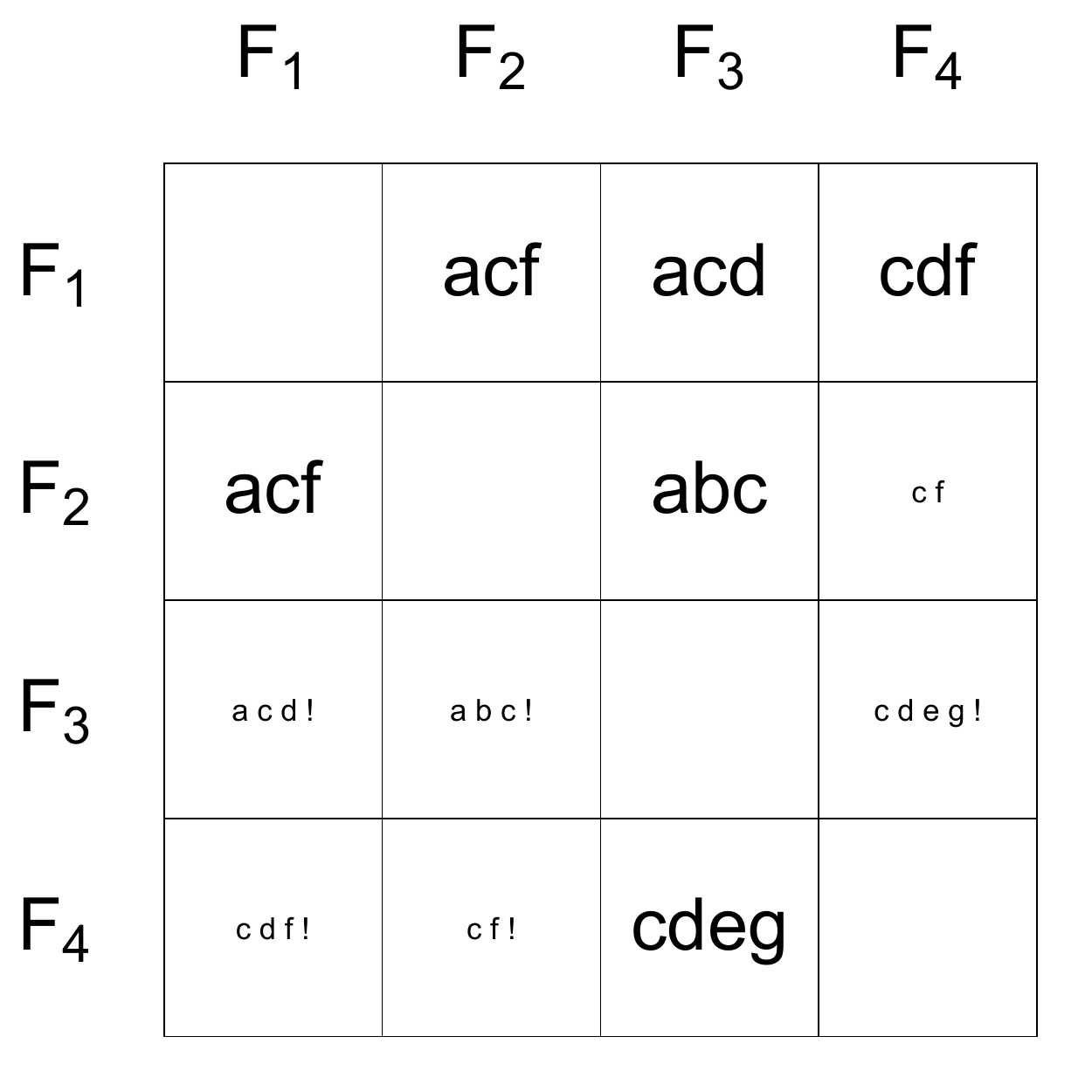}

{\it Table 2: All facet intersections (including cops and hooligans) for the complex $\Delta_0$}

If $\Delta_0$ has a shelling $(G_1,G_2,G_3,G_4)$ then $F_3=G_1$ because $F_k\cap F_3\not\prec F_3$ (recall condition (1)). Hence any potential shelling is of type
$(F_3,G_2,G_3,G_4)$. Neither $F_1$ nor $F_2$ can preceed $F_4$ because $F_1\cap F_4\not\prec F_4$ and $F_2\cap F_4\not\prec F_4$ and $F_i\cap F_4\not \subseteq F_3\cap F_4\ (i=1,2)$. But $F_3$ may preceed $F_4$ since
$F_3\cap F_4\prec F_4$. Hence any potential shelling is of type $(F_3,F_4,G_3,G_4)$. If $G_3=F_2$ then $F_4\cap F_2=\{c,f\}\not\prec F_2$, and the only other predecessor  of $F_2$ in $(F_3,F_4,F_2,F_1)$ is $F_3$ which doesn't save the situation since $F_4\cap F_2=\{c,f\}\not\subseteq\{a,b,c\}=F_3\cap F_2$. Our last hope remains
$(F_3,F_4,F_1,F_2)$. As to $F_1$, no position of $F_1$ could ever cause harm because of $F_k\cap F_1\prec F_1$ for all $k\neq 1$. As to $F_2$, the intersection
$F_4\cap F_2=\{c,f\}$ has now {\it ceased} to be cumbersome because now $F_4\cap F_2=\{c,f\}\subseteq\{a,c,f\}=F_1\cap F_2\prec F_2$. To summarize, $(F_3,F_4,F_1,F_2)$ is the unique shelling of $\Delta_0$. We shall often only write the indices of a shelling, i.e. write $(3,4,1,2)$ or even 3412 for the shelling above.

{\bf 3.3} Let $F_1,\ldots,F_n$ be the facets of some complex $\Delta$. If $A\subseteq [n]$ and $k\in [n]$ we shall ponder this condition:

(10)\quad $(\forall i\in A)\ (\exists j\in A)\ F_i\cap F_k\subseteq F_j\cap F_k\prec F_k$

Suppose $(i_1,\ldots,i_s,k)$ is a partial shelling\footnote{Recall that this is shorthand for $(F_{i_1},\ldots,F_{i_{s}},F_k)$.} of $\Delta$, and $A:=\{i_1,\ldots,i_{s}\}$. Then $(A,k)$, in view of (1), satisfies (10). Thus if ${\cal P}$ is the family of all pairs $(A,k)$ satisfying (10), then ${\cal P}$ satisfies (4a) and (4b). A moment's thought shows that (4c) holds as well. Hence all these $(A,k)$ are potential setment-suffices,  and ${\cal P}$ is a PSS-poset for the language ${\cal L}\subseteq [n]^*$ of all partial shellings of $\Delta$.

Here comes the systematic procedure to calculate this PSS-poset for given facets $F_1,\ldots,F_n$  of a complex $\Delta$. The suffix $k\in [n]$ being fixed, focus on  the family of all PSSes of type $(A,k)$. In view of (1) call $j\neq k$ a {\it cop} (for $k$) if  $F_j\cap F_k\prec F_k$, and call $i\neq k$ a {\it hooligan} (for $k$) if $F_i\cap F_k\not\prec F_k$.
 As illustrated in our toy example, the potential setment $A$ of $k$ may contain a hooligan $i$ only if $A$ contains at least one cop $j$ {\it policing} $i$ in the sense that $F_i\cap F_k\subseteq F_j\cap F_k\prec F_k$. In this situation we also  say: $i$ {\it has} cop $j$. We distinguish two types of hooligans; the {\it policeable} ones (i.e. having at least one cop), and  the {\it non-policeable} ones. For instance, if $k$ has no cops, then all hooligans are non-policeable.
For each fixed $k$ the following holds:

(11)\quad There are cops $j$ for $k\ \Leftrightarrow\ $ There are nonempty potential setments $A$ of $k$

To verify (11), let $j$ be a cop for $k$. Then $A:=\{j\}$ is a potential setment of $k$. Conversely, when $A$ is any potential setment of $k$, we claim there is a cop in $A$. Indeed, otherwise (since $A\neq\emptyset$ by assumption) there would be an $i\in A$ violating (1). This proves (11).

{\bf 3.4} Suppose some $n=9$ facets $F_i$ are such that the potential setments $A$ of $k=9$ are subject to these (and only these) conditions:

\begin{description}
\item[(12a)] the hooligan 8 has cop 5;
\item[(12b)] the hooligan 4 has cops 3, 5, 7;
\item[(12c)] the hooligans 1 and 2 have no cops;
\item[(12d)] the cop 6 is not policing any hooligans.
\end{description}

 Identifying subsets of $E=[8]$ with length 8 bitstrings the required potential setments $A$ are exactly the models of $\varphi$ in (7), and hence can be calculated as in 3.1. Besides compression, which is the raison-d'\^{e}tre of the $e$-formalism, in our scenario the $e$-algorithm benefits from two pleasant circumstances. 

First, it is {\it independently} run $n$ times, i.e. once for  each suffix $k\in [n]$. (Recall that the three $012e$-rows in Table 1 result from the {\it particular} value $k=9$.) For each individual $k$ the task is benign. Specifically, for each $k\in [n]$ its dual Horn formula has  as many clauses as policeable hooligans, i.e. at most $n$. The fewer cops, the more nonpoliceable hooligans one tends to have. 

Second, the structure of our clauses trigger fewer $012e$-rows than an equal number of random clauses would. This is because all positive literals derive from one subset (the cops) and all negative literals from a {\it disjoint} one (the nonpoliceable hooligans).

{\bf 3.5} It follows from (11) that when $k\in [n]$ has no cops then $k$ must be the first element in any shelling of $\Delta$. Consequently, if more than one $k\in [n]$ have no cops, then $\Delta$ has no shellings. In this case we speak of a {\it Type 1 failure of shellability}. 
 Suppose $\Delta$ is such that each $k\in [n]$ has at least one non-policeable hooligan for $k$. Then no $k$ can be the last element in a potential shelling, i.e. $\Delta$ is unshellable. We call this a {\it Type 2 failure of shellability}. 
Type 1 and 2 are not based\footnote{But they carry over
as follows to general PSS-posets. Type 1 occurs iff there are suffixes $k\neq k'$ such that PSSes $(A,k)$ and $(A,k')$ exist at most when $A=\emptyset$. Type 2 occurs iff for each $k\in [n]$ all PSSes $(A,k)$ have $|A|\leq n-2$. } on the PSS concept and
 are visualized in Table 2. Thus in the row of $F_k$ an intersection $F_j\cap F_k$ is written large if $j$ is a cop (for $k$), and small if $j$ is a hooligan. Furthermore, each nonpoliceable hooligan is accompanied by '!'.
Accordingly a Type 1 failure occurs iff at least two rows consist entirely of small entries, and a Type 2 failure occurs iff each row features some '!'. It is possible that both types of failure occur simultaneously. In contrast two further types of failure are naturally described in terms of PSSes. Namely we speak of a {\it Type 3 failure of full words} if there are no PSSes $(A,k)$ with $|A|=n-1$. Furthermore, a {\it Type 4 failure of full words} is simply the absence of
admissible $n$-chains. Recall from Proposition 1 that this is sufficient {\it and} necessary for non-existence of full words.
 
{\bf 3.6} As is well known, depth-first-search (DFS) computations invite distributed computing. For instance, instead of having one initial stack (A) in Figure 2, one could distribute its ingredients to three stacks, process these stacks independently, and then collect the results. Being based on DFS (as detailed in [W2]) the $e$-algorithm of 3.1 also lends itself to distributed computing. In particular observe that only the $e$-algorithm, and no other fuzz, is required to detect a Type 3 failure of full words.

\section{Matroid complexes, chessboard complexes, and more}

We first stick with shellability but specialize to matroid complexes (4.1) and the related chessboard complexes (4.2). Afterwards, in 4.3, we turn to peelings of arbitrary set systems. We chose the term 'peelings' to indicate a superficial proximity to shellings. Although calculating  PSS-posets for peelings works along different lines, the $e$-algorithm continues to be useful. In 4.4 peeling languages are seen to generalize poset languages (whose basic words are the linear extensions of the poset). Other uses of PSS-posets, e.g. for so-called alternative precedence languages,  are sketched.

{\bf 4.1} A {\it matroid complex} is the simplicial complex $\Delta\subseteq 2^W$ of all independent sets of a matroid on the  ground set $W$. It is well known that each permutation of $W$ induces a special type of shelling of $\Delta$. However, as mentioned in the introduction the number of shellings $ns(\Delta)$ can be much larger than
 $(|W|!)\ $. Determining  $ns(\Delta)$ gets interesting already for rank 2 matroids; let us focus on the unique {\it simple} rank 2 matroid $M(2,m)$ with groundset $[m]$. Then the $n={m\choose 2}$   facets of the coupled matroid complex $\Delta(M(2,m))$ are all 2-sets of $[m]$, and so the shellings of $\Delta(M(2,m))$ bijectively match those listings of the edges of the complete graph $K_m$ for which all $n$ 'initial' induced subgraphs are connected. For instance, if $m=5$ then all 7! many length $n$ listings beginning with $\{1,2\},\ \{2,3\},\ \{3,4\}$ {\it are} shellings, whereas all listings beginning with $\{1,2\},\ \{2,3\},\ \{4,5\}$ are {\it no} shellings. Straightforward reasoning along these lines shows that $ns(\Delta(M(2,m)))\, /\, n!$
equals $1,\ 4/5,\ 4/7$ for $m=3,4,5$, and that this proportion goes to $0$ as $m$ goes to infinity. Modulo the probabilities $p(s,t)$ that a random graph with $s$ vertices and $t$ edges is connected, it is an exercise to give an explicite formula for $ns(\Delta(M(2,m)))$. It is a bit clumsy, not least because only recursive formulas for $p(s,t)$ seem to be known. We thus calculated $ns(\Delta(M(2,6)))=498161664000$ brute-force with our algorithm (Table 3).

\begin{tabular}{|c|c|c|c|c|c|c|} 
complex & $\#\,$facets & $\#\,$PSSes & $\#\,012e$-rows & Time & $\#$ shellings\\   \hline
 & & & & & & \\ \hline
 $\Delta(M(2,6))$ & 15 & 244 800 & 15 & 16387  & 498 161 664 000\\
$\Delta(PM(3,3,2))$ & 18 & 1 884 672 & 90 & 10064  & 14 004 606 481 920\\
$\Delta(PM(2,2,2,2))$ & 16 & 270 336 & 208  & 85 & 6 163 021 824\\
$\Delta(M(K_4))$ & 16 & 470400 & 80 & 8308  & 722 965 625 856\\
$\Delta(CB(4,4,4))$ & 24 & 110 100 480 & 120 & 405  & 0\\
$\Delta(CB(3,2,2,2,1))$ & 16 & 419 328 & 34 & 7155 &  194 527 872 000\\
$\Delta(CB(3,3,2,2,1))$ & 20 & 7 274 496 & 72 & 32463  & 116 916 202 200 752\\
$\Delta(CB(4,3,2,1))$ & 14 & 59904 & 43 & 21  & 44 176 168\\    \hline

\end{tabular}

{\it Table 3: The number of shellings for various simplicial complexes}

Before we continue with matroids, let us comment on the structure of Table 3. The first column identifies the simplicial complex according to the notation introduced in Section 4. Column 2 gives the number $n$ of its facets. The larger $n$ the larger the number of PSSes (column 3) tends to be. More important than the sheer size of the PSS-poset is how good it got compressed with $012e$-rows (column 4). For instance, the 20 facets of $\Delta=\Delta(CB(3,3,2,2,1))$ induced a large PSS-poset but only\footnote{Running the $e$-algorithm took a mere 0.16 sec and $<0.3$ sec for all instances in Table 3.} seventy two $012e$-rows, each one of which accomodating an average of 100 000 PSSes.  Thus when 
all good upper covers of $(A,k)$  must be found (recall, they have setment $A\cup\{k\}$) we need not scan 7 million candidate PSSes but only 72 multivalued rows. Deciding whether $ns(\Delta)=0$ or $ns(\Delta)>0$ is faster\footnote{Observe that on the 'way up' in the PSS-poset (see 2.4) the  calculation of $ns(\Delta)$ demands to {\it additionally} store the good upper covers of each good PSS. This led to some clumsy (but perhaps amendable) book keeping. The 'way down', i.e. the adding up of numbers, is fast: Only 13 of the 32463 seconds were spent that way.} than obtaining the exact count $ns(\Delta)=116916202200$ (column 6) which took 32463 seconds (column 5).

As to higher rank  matroids, a natural class to investigate are partition matroids. For instance we denote by $PM(\{1,3,7\},\{2,5,8\},\{4,6\})$ the the partition matroid on $W=[8]$ whose bases are the twelve 3-element transversals of the sets $\{1,3,7\},\ \{2,5,8\},\ \{4,6\}$. Since often only the cardinalities of these sets matter, we can alternatively
 write $PM(3,3,2)$ for this matroid, and $\Delta(PM(3,3,2))$ for its matroid complex. Apart from this kind, Table 3 also features the number of shellings of the matroid complex induced by the rank 3 polygon-matroid $M(K_4)$ on the edge set of the complete graph $K_4$.

 {\bf 4.2} The intersection $\Delta_1\cap\Delta_2$ of two matroid complexes $\Delta_1$ and $\Delta_2$ needs not be a matroid complex, and so the shellability of $\Delta_1\cap\Delta_2$ is not guaranteed. The special case where $\Delta_1$ and $\Delta_2$ derive from partition matroids has attracted particular attention. Specifically, let $M_1=PM(A_1,\ldots,A_s)$ and $M_2=PM(B_1,\ldots,B_t)$ be based on $W$ and such that $|A_i\cap B_j|\le 1$ for all $i\in [s]$ and $j\in [t]$. Then, as is well known, the elements of $W$ can be identified with the squares of a (possibly perforated) 'chessboard' in such a way that the $A_i$'s match the rows, and the $B_j$'s match the columns. Furthermore the faces of $\Delta_1\cap\Delta_2$ (where $\Delta_i=\Delta(M_i)$) bijectively match the configurations of non-taking rooks on the chessboard. For instance, if $\Delta_1=\Delta(PM(\{1,2,3,4\},\{5,6\},\{7,8\},\{9,10\}))$ and $\Delta_2=\Delta(PM(\{1,6,8\},\{2,5,9\},\{3,7\},\{4,10\}))$ then the ensuing chessboard is rendered in Figure 5(a). The {\it chessboard complex} $\Delta=\Delta_1\cap \Delta_2$  is not pure, i.e. 12 facets are not equicardinal (see Figures 5(b),(c), and similarly on (d)).

In the sequel we look at more regular (unperforated) chessboards $CB(i_1,i_2,\ldots,i_s)$ and their complexes $\Delta(CB(i_1,i_2,\ldots,i_s))$. Namely, by definition $CB(i_1,i_2,\ldots,i_s)$ is a gluing of contiguous left-aligned rows of of lengths $i_1\ge i_2\ge \ldots\ge i_s$. Figure 5(d) displays $CB(4,4,2,2)$. 
In [Z] it is shown that chessboard complexes of type $\Delta(CB(k,k,\ldots,k))$ are (vertex-decomposable and whence) shellable whenever $k$ is large enough. Specifically, if there are $j$ rows of length $k$ then $k\ge 2j-1$ guarantees shellability. This condition is e.g. satisfied for $\Delta(CB(5,5,5))$ but not for $\Delta(CB(4,4,4))$. Indeed, it turns out that $ns(\Delta(CB(4,4,4)))=0$, the largest partial shellings having cardinality 13 instead of 24.
In  [CZ] we learn that complexes of chessboards with at least $m$ rows of length $m$, and $m-1$ rows of length $m-1$,..., and at least 1 row of length 1, are (vertex-decomposable and whence) shellable. We establish that otherwise the complex may still be shellable, such as $\Delta(CB(3,3,2,2,1))$ and $\Delta(CB(3,2,2,2,1))$ in Table 3. On the other hand, the chessboard complexes with row patterns (3,3,2,1,1) and (3,3,2,2) narrowly fail to be shellable, and a blunt Type 1 failure (4 facets without cops) occurs for (4,4,2,2). While so called 
{\it Stirling complexes} $\Delta(CB(m,m-1,\ldots,2,1))$ drastically violate the above sufficiency criterion, it turns out that $\Delta(CB(4,3,2,1))$ is nevertheless
shellable.

\hspace*{-1.3cm}\includegraphics[scale=0.62]{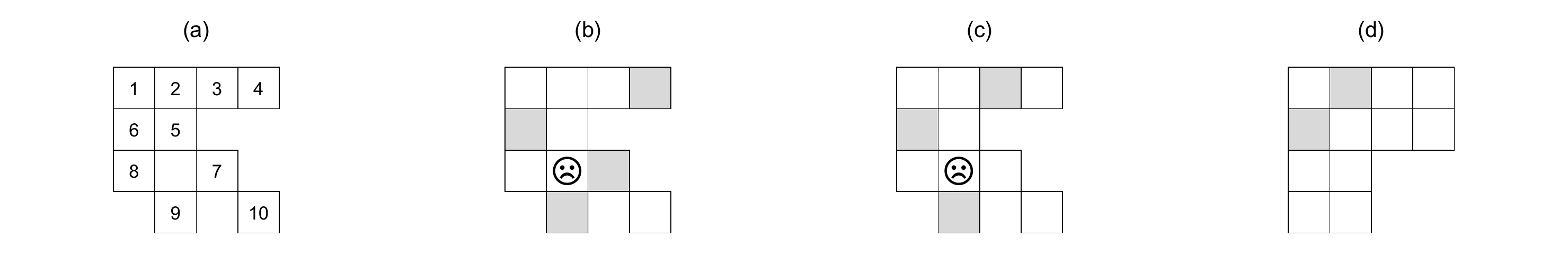}

{\it Figure 5: Two chessboards with some maximal non-taking-rook-configurations}

{\bf 4.3} Some preliminaries are in order before we get to peelings. Given a target set $X$ and a family ${\cal G}$ of subsets, a {\it set covering} (of $X$) is a subfamily ${\cal F}\subseteq {\cal G}$ such that $X \subseteq 
\bigcup {\cal F}$. For each $x \in X$ put ${\cal G}_x: = \{Y \in {\cal G}: x \in Y\}$. As is obvious and well known, ${\cal F}$ is a set covering iff ${\cal F}$ is a transversal of the set system $\{{\cal G}_x: x \in X\}$. The $e$-formalism of 3.1 can also be exploited [W2] to achieve a compressed representation of the family of transversals of $\{{\cal G}_x: x \in X\}$.

Let now $F_1,\ldots,F_n$ be any sets. Despite notation we do not assume that they are the facets of a simplicial complex, i.e. they need not be incomparable. A permutation $(G_1,\ldots,G_n)$ of $F_1,\ldots,F_n$ is called\footnote{This is not standard terminology. We chose 'peeling' because it conveys a proximity to 'shelling'. }  a {\it peeling} if for all $2\leq k\leq n$ it holds that

(13)\quad $(G_1\cup\cdots\cup G_{k-1})\cap G_k\prec G_k$

Albeit (13) may look like a falsely remembered version of shellability (1), one can motivate studying this condition. If say $n=4$ and $(G_1,G_2,G_3,G_4)=(F_2,F_3,F_4,F_1)$ satisfies (13) then 2341 is a basic word of the {\it peeling language} ${\cal L=L}(F_1,..,F_4)$ on the alphabeth $E=[4]$. In order to calculate a PSS-poset ${\cal P(L)}$ of the peeling language ${\cal L}$ induced by $F_1,\ldots,F_n$ we proceed as follows. 

Similarly to shellings we reduce the task to the calculation of the families $PSSes[k]$ of all PSSes with suffix $k$. Because of (13) for any such $(A,k)$ in $PSSes[k]$ it holds that

(14)\quad $(\bigcup \{F_i:\ i\in A\}) \cap F_k\prec F_k$.

 Obviously $PSSes[k]$ is the disjoint union of the families $PSSes[k,v]$ where $v$ ranges over $F_k$, and where $PSSes[k,v]$ consists of all $(A,k)$ with

(15)\quad $(\bigcup \{F_i:\ i\in A\}) \cap F_k= F_k\setminus \{v\}$.

Putting $X: = F_k \setminus \{v\}$ and letting ${\cal G}$ be the family of all sets $F_i$  with $v \not\in F_i$ we see that each set covering ${\cal F} \subseteq {\cal G}$ of $X$ induces an $(A,k)\in PSSes(k,v)$ that has $A:=\{i:\ F_i\in {\cal F}\}$. Conversely, every member of $PSSes(k,v)$ arises in this way. Recalling the remarks about set coverings one sees how a compressed representation of $PSSes[k,v]$ by means of $012e$-rows can be achieved. Collecting the rows for all $k\in [n]$ and $v\in F_k$ yields the desired compressed representation of ${\cal P(L)}$.

{\bf 4.4} Let $({\cal P},<)$  be a poset and let $\{F_1,\ldots,F_n\}$ be the set system of all its principal ideals $p\!\!\downarrow :=\{b\in P:\ b\leq p\}$. Let $(p_1\!\!\downarrow,\ldots,p_n\!\!\downarrow)$ be any peeling of $\{F_1,\ldots,F_n\}$. Thus say $(p_1\!\!\downarrow \cup\, p_2\!\!\downarrow)\cap p_3\!\!\downarrow=
p_3\!\!\downarrow\setminus\{b\}$. Suppose we had $b\neq p_3$. Then $p_3\in p_1\!\!\downarrow \cup\, p_2\!\!\downarrow$, say $p_3\in p_2\!\!\downarrow$. This implies $b<p_3<p_2$, and thus the contradiction  $b\in (p_1\!\!\downarrow \cup\, p_2\!\!\downarrow)\cap p_3\!\!\downarrow$. Therefore $(p_1\!\!\downarrow \cup\, p_2\!\!\downarrow)\cap p_3\!\!\downarrow=p_3\!\!\downarrow\setminus\{p_3\}$.
 Similarly $(p_1\!\!\downarrow \cup\ldots\cup p_{k-1}\!\!\downarrow)\cap p_k\!\!\downarrow=p_k\!\!\downarrow\setminus\{p_k\}$. This amounts to say that $(p_1\!\!\downarrow,\ldots,p_n\!\!\downarrow)$ is a linear extension of $({\cal P},\leq)$. Conversely, every linear extension yields a peeling\footnote{But not every peeling language arises from a suitably chosen poset. For instance not ${\cal L}={\cal L}_1$ which is induced by $F_1=\{a,b,e\},\ F_2=\{d,f\},\ F_3=\{b,c,e\},\ F_4=\{d,e\}$.} of $\{F_1,\ldots,F_n\}$. We denote this type of peeling language as ${\cal L(P)}$.

Each setment of a partial linear extension is an order ideal, and it is well known that conversely every order ideal arises this way. Thus we can skip `potential' setment-suffixes and focus on setments (=order ideals) right away. This has been done in [W3]. 
We mention that  ${\cal L(P)}$ is a greedoid language in the sense of 2.1. In fact, it is an example of a wider class of greedoid languages, so-called alternative precedence languages. It turns out that similar to order ideals the setment-suffixes of alternative precedence languages can be targeted directly (work in progress).

 Many other ordering problems, e.g. concerning integers $a_1,\ldots,a_n$ subject to arithmetic restrictions, can be tackled by the use of PSS-posets. For instance, the reqirement could be that for each $k\in [n]$ the sum $a_1+\cdots+a_{k-1}$ is congruent to $a_k$ modulo $7$.

 \section*{References}
 
\begin{enumerate}
\item [{[B]}] A. Bj\"{o}rner, The homology and shellability of matroids and geometric lattices, chapter 7 in Matroid Applications, Cambridge University Press 1992.
\item [{[BZ]}] A. Bj\"{o}rner, G. Ziegler, Introduction to greedoids, chapter 8 in Matroid Applications, Cambridge University Press 1992.
\item [{[CZ]}] E. Clark, M. Zeckner, Simplicial complexes of triangular Ferrers boards, J. Algebraic Combin. 38 (2013), no. 1, 1–14. 
\item [{[MH]}] S. Moriyama, M. Hachimori, h-assignments of simplicial complexes and reverse search, Disc. Appl. Math 154 (2006) 594-597.
\item[{[W1]}] M. Wild, Compactly generating all satisfying truth assignments of a Horn formula, Journal on Satisfiability, Boolean Modeling and Computation 8 (2012) 63-82.
	\item [{[W2]}] M. Wild, Counting or producing all fixed cardinality transversals. Algorithmica (2014) 69: 117-129.
	\item[{[W3]}] M. Wild, An efficient data structure for counting all linear extensions of a poset, calculating its jump number and the likes,
	               preliminary version on ResearchGate.
	\item[{[W4]}] M. Wild, ALLSAT compressed with wildcards: Partitionings and face-numbers of simplicial complexes, submitted. 							
	\item[{[Z]}] G. Ziegler, Shellability of chessboard complexes, Israel J. Math. 87 (1994) 97-110.
\end{enumerate}

\end{document}